\documentclass[11pt]{article}
\usepackage{stmaryrd}
\usepackage{amssymb}
\usepackage{color, amsmath,amssymb, amsfonts, amstext,amsthm, latexsym}

\usepackage{amssymb, epsfig, amssymb, latexsym}
\usepackage{amsmath}
\usepackage{graphicx}
\usepackage{longtable}
\usepackage{graphicx}
\usepackage{subfigure}
\usepackage{booktabs}
\usepackage{cite}

\usepackage[colorlinks=false, linktocpage=true]{hyperref}

\renewcommand{\phi}{\varphi}

\begin{document}

\title{Quantifying model uncertainty for  the observed non-Gaussian data by the Hellinger distance }

\author {Yayun Zheng$^{1}$, Fang Yang $^{*2}$, Jinqiao Duan$^{3}$,  J\"urgen Kurths$^{4, 5}$\\
\\$^1$ School of Mathematical Sciences, Jiangsu University, Zhenjiang 212013, China.
\\$^2$ School of Mathematics and Statistics and Center for Mathematical Science, \\ Huazhong University of Science and Technology, Wuhan, 430074, China.
\\$^3$ Department of Applied Mathematics, Illinois Institute of Technology,\\  Chicago, IL 60616, USA.
\\$^4$ Potsdam Institute for Climate Impact Research, Potsdam, 14473, Germany. 
\\ $^5$Department of Physics, Humboldt University, Berlin, 12489, Germany.
\\$^6$ Centre for Analysis of Complex Systems, Sechenov First Moscow State Medical University, \\ Moscow, Russia. 
}
 \footnotetext[1]{* Corresponding author.\\
 E-mail address: yfang@hust.edu.cn (Fang Yang)} 

\date{}

\maketitle

\begin{abstract}
Mathematical models for complex systems under random fluctuations often certain uncertain parameters. 
However, quantifying model uncertainty for  a  stochastic differential equation  with an $\alpha$-stable L\'evy process is still  lacking. Here,  we propose an approach  to infer  all the uncertain  non-Gaussian parameters  and other system parameters by minimizing the Hellinger distance over the parameter space. The Hellinger distance measures  the similarity between an empirical probability density of non-Gaussian observations and a solution (as a probability density) of the associated nonlocal Fokker-Planck equation.
Numerical experiments verify that our method is feasible for estimating single and multiple parameters. Meanwhile, we  find an optimal estimation interval of the estimated parameters.  This  method is beneficial for extracting governing  dynamical system models under non-Gaussian fluctuations, as in the study of abrupt climate changes in the Dansgaard-Oeschger events. 


\end{abstract}

\section{Introduction}

Complex systems under influences of  random fluctuations also have  uncertain parameters \cite{bishwal2007parameter}.  An important problem in modeling such random processes by stochastic differential equations (SDEs) is to estimate uncertain parameters from observations of the stochastic paths. 


A Brownian motion has properties of continuous sample paths, normal diffusion and light tail (probability density  decays exponentially),  theoretical results on 
parametric estimations  for SDEs driven by Brownian motion are relatively well developed.
The Gaussian kernel density estimator \cite{iacus2009simulation, botev2010kernel} and the Bayesian estimator\cite{papaspiliopoulos2012, beaumont2002approximate} are well-known approaches  for parameters estimation of a drift function when the observations of  complete paths are available.  
Furthermore, by the  nonparametric estimation method  of  Kramers-Moyal coefficients \cite{honisch2011publisher}, the statistical definitions of conditional  first and second moments  \cite{siegert1998analysis}  or variational formulation of the stationary Fokker-Planck equation \cite{batz2016variational}, we could provide an expression for the drift function and the diffusion one.

However, various complex phenomena  involve non-Gaussian fluctuations, with properties such as  intermittent jumps, anomalous diffusion, and heavy tail (probability density  decays with power law) distribution. A heavy-tailed distribution, like a L{\'e}vy distribution, is characterized by a high likelihood for extreme events, compared to a normal distribution. For instance, Ditlevsen shows that the  paleoclimatic records for Dansgaard-Oeschger events have a strong non-Gaussian distribution \cite{ditlevsen1999}.  The protein production occurs in bursts which are observed during a genetic regulation \cite{cai2006stochastic}. 
Meanwhile,  experimental studies find that  L{\'e}vy flights are  an optimal pattern when the prey is sparsely and randomly distributed for open-ocean predatory fish\cite{humphries2010environmental}.
The L{\'e}vy process  is also used in other scientific domains, for example, it has been shown that certain stock price  has statistical properties that are compatible with a L{\'e}vy random walk\cite{mantegna1991levy}. Additionally, in the field of cognitive research, a few studies provide evidence of L{\'e}vy  processes, e.g., to  search and cluster in semantic memory \cite{montez2015role} and human decision making\cite{voss2019sequential}.  \textcolor[rgb]{0,0,0}{An $\alpha$-stable L\'evy process is thought to be an appropriate model for a non-Gaussian heavy-tailed process. 
Therefore, for modeling these complex systems, it becomes  necessary and significant to consider  parameters estimation for a stochastic system driven by the $\alpha$-stable L\'evy process. In general, the $p$-th moment of an $\alpha$-stable  L\'evy random variable is finite if and only if $p < \alpha$ ($ 0 < \alpha <2$), so it does not have second moments. Meanwhile, the stationary probability density of an $\alpha$-stable L\'evy process does not always exist. Due to these disadvantages, unfortunately, the existing methods for  parameter estimation  of Brownian motion are not applied for dynamical system with non-Gaussian fluctuations. }

%


 \textcolor[rgb]{0,0,0}{There are few results about the parametric estimation  for stochastic processes driven by  L\'evy processes. 
In some special cases, it is possible to  infer parameters only for the drift function assuming that the values of other parameters are known.
For example,  a simple Ornstein-Uhlenbeck process is considered,  i.e., the drift function is known to be linear,  or stochastic processes are driven by a compound Poisson process \cite{ogihara2011quasi}. In these works, the quasi-maximum likelihood, a self-weighted least absolute deviation estimator \cite{masuda2010approximate} or trajectory fitting estimator are established for discretely observed L\'evy  processes.
For an $\alpha$-stable L\'evy process,  the parametric estimation problem becomes more difficult because 
the second moment does not exist. Recently, Hu and Long \emph{et al.} \cite{hu2007parameter, hu2009least} addresses a trajectory fitting  and a least-square estimator on estimation of a drift parameter for a stochastic system under an $\alpha$-stable L\'evy noise.  Fasen \cite{fasen2013statistical} extended the results to  high dimensions.}


In the above-mentioned works, one can only estimate the drift parameters. Meanwhile, the non-Gaussian index $\alpha$  plays a decisive role in the construction of  L\'evy processes. 
An $\alpha$-stable L\'evy process has larger jumps with lower jump probabilities when $\alpha$  is small ($ 0 < \alpha <1$), while it has smaller jumps with higher jump frequencies for large $\alpha$
 values  ($1 < \alpha <2$). 
 The special cases for $\alpha  = 1$ and $\alpha  = 2$ correspond to  the Cauchy process and the Brownian motion, respectively.  Therefore, the estimation of the parameter $\alpha$ is extremely important. 
 There are some  simple and straightforward approaches to learn this $\alpha$ from the path observation, such as the slope of the log-log linear regression \cite{gopikrishnan1999scaling} or the Hill estimator \cite{hill1975simple}. These methods do not assume a parametric form for the entire distribution function, but focus only on the tail behavior. However,  the true tail behavior of  L{\'e}vy distribution is visible only for extremely large data sets, or it is a challenge to choose the right value of the largest order statistics.
 
There have been no available estimators simultaneously for the drift parameter and other $\alpha$-stable L\'evy parameters, including $\alpha$ and other non-Gaussian parameters. An alternative method is relied on the characteristic function of $\alpha$-stable L\'evy process \cite{yang2011}. Based on the ergodic theory and sample characteristic functions, Cheng \emph{et al.} \cite{cheng2020generalized} study a Ornstein-Uhlenbeck process with the $\alpha$-stable L\'evy noise. The parameter estimation for  $\alpha$ and the other parameters is obtained by matching the empirical characteristic function with the corresponding theoretical one. We note that  a method of numerical optimization is devised in \cite{gao2016quantifying}, where  two deterministic quantities: mean exit time or the escape probability is observed to estimate the uncertain parameter and other system parameters. It is based on solving an inverse problem for a deterministic, nonlocal partial differential equation.

\textcolor[rgb]{0,0,0}{The existing works provide certain approaches to estimate the drift parameter and other $\alpha$-stable L\'evy parameters. With severe limitations,  (i) the drift term can only be  a linear function; (ii) the empirical characteristic function is approximately defined; (iii) it is  difficulty to observe mean exit time or first escape probability from the discrete time series data. In response to the existing challenge, we are interested in finding an effective and feasible approach for parameter estimation of a stochastic system under an $\alpha$-stable   L\'evy noise. The method can be applied for a nonlinear drift term, $\alpha$ and other system parameters  can be estimated simultaneously. Compared with the characteristic function and other quantities, a probability density or probability distribution becomes easy to obtain  from an observation data set with heavy-tailed distribution. }

\textcolor[rgb]{0,0,0}{We recall some recent works on estimating parameters of $\alpha$-stable  L\'evy  processes based on probability densities. 
Chen and Chen \cite{chen2020maximum} choose a mixture of Cauchy and Gaussian distribution to approximate the probability density function of the $\alpha$-stable  Ornstein-Uhlenbeck distribution. By means of transition function and Laplace transform, they construct an explicit approximate sequence of the maximum likelihood function to obtain the estimation of parameters.  Inspired by the derivation of the differential Chapman-Kolmogorov equation, Li and Duan \cite{liyangdata} derive Kramer-Moyal formulas to express the jump measure, drift and diffusion coefficient of a stochastic differential equation with respect to the transition probability density .}

The probability density function  for SDE driven by an $\alpha$-stable L\'evy process satisfies a deterministic, nonlocal differential equation with an initial condition, i.e., nonlocal Fokker-Planck equation, which have a nonlocal or fractional Laplacian term. 
 In terms of theory, we derived the Fokker-Planck equations for   Marcus SDEs driven by L\'evy processes in high dimensional\cite{sun2017marcus}. In terms of numerical calculations, taking advantage of the Toeplitz matrix structure of the time-space discretization,  Gao \emph{et al.} \cite{gao2016fokker} proposed a fast and  accurate numerical algorithm   to simulate  nonlocal Fokker-Planck equations under either absorbing or natural conditions. \textcolor[rgb]{0,0,0}{Meanwhile, a piecewise integro quadratic spline interpolation approach\cite{moghaddam2019, keshi2019numerical, moghaddam92integro} and a finite element method  \cite{nie2020numerical} are developed for  the approximate nonlocal or fractional  integral. }
 
\textcolor[rgb]{0,0,0}{ Consequently,  we propose an approach  to infer  simultaneously  the drift parameter and other $\alpha$-stable L\'evy parameters.
Our method is based on minimizing the Hellinger distance between the observed probability distribution and the solution (as a probability distribution) of the associated Fokker-Planck equation for a general stochastic dynamical system driven by  an $\alpha$-stable L\'evy process.}

In the present paper,  we consider the parameter estimation problem of an $\alpha$-stable   L\'evy stochastic dynamical system containing uncertain parameters.
 In Section 2, we propose a method of estimating the uncertain parameters based on the Hellinger distance of the probability densities. In Section 3,  we  present some simulation results of  estimation for  single and multiple parameters by minimizing the Hellinger distance.
Finally we give some concluding and future works  in Section 4.


\section{Methods}
We consider  a  dynamical system with heavy-tailed uncertainty, which could be modeled by a stochastic process $X(t) $
\begin{equation}\label{p0}
{\rm d}X(t) = f(X(t), \theta){\rm d}t+ \epsilon {\rm d}L^\alpha(t) , \quad X(0)=x_0 \in \mathbb{R}^1,
\end{equation}
where the drift function $f(x, \theta)$ has the uncertain parameter $\theta$,   and a scalar  symmetric $\alpha$-stable L\'{e}vy process $L^\alpha(t)$ with the non-Gaussian index $0 < \alpha< 2$ is defined in a probability space $(\Omega, \mathcal{F}, \mathbb{P})$. The parameter $\epsilon$ is the non-negative $\alpha$-stable L\'{e}vy noise intensity.

A scalar symmetric $\alpha$-stable L\'{e}vy process  is characterized by a generating triplet ($b, Q, \nu_\alpha$),  a linear coefficient $b$, a diffusion parameter  $Q$, and a nonnegative Borel measure $\nu_\alpha$. This  jump measure  $\nu_\alpha$ is defined on $(\mathbb{R}^{1}, \mathfrak{B}(\mathbb{R}^{1}))$\cite{applebaum2009levy} by:
\begin{equation*}
 \nu_{\alpha} = \frac{C_\alpha dy}{|y|^{1 + \alpha}},
\end{equation*}
with $0 < \alpha < 2$ and $C_\alpha = \frac{\alpha}{2^{1 - \alpha}\sqrt{\pi}}\frac{\Gamma(\frac{1+\alpha}{2})}{\Gamma(1-\frac{\alpha}{2})}$.  In this paper, we consider an $\alpha$-stable L\'{e}vy process with a triplet $(0, 0, \nu)$, i.e., a pure jump process.

 For $0<\alpha<2$,  the $\alpha$-stable L\'{e}vy process $L^\alpha(t)$ has a heavy-tailed  distribution\cite{sato1999levy}
\begin{equation}
\mathbb{P}(|L^\alpha(t)|>y) \sim \frac{1}{y^\alpha}, \notag
\end{equation}
as the tail estimate decays in a power law. Therefore $\alpha$ is also called the power parameter.  The tail behavior is different from the Brownian motion with light tail, as the tail  decays exponentially.

We assume that the drift term $f$ is  local Lipschitz continuous. Then the SDE (\ref{p0})  has  a unique  solution \cite{applebaum2009levy}. The conditional probability density $p(x, t | x_0, 0) \triangleq p(X(t) = x | X(0) = x_0)$ represents the density of the $X(t)$  given a value $x_0$ at initial time.
 For convenience, we drop the initial condition and simply denote it by $p(x, t)$. There exists sufficient condition for the existence and regularity of the probability density $p(x, t)$ for some SDEs driven by L\'{e}vy processes. The existence is based on Malliavin calculus with jumps under H\"ormander's condition, see Refs. \cite{zhang2014densities, song2015regularity, chen2015stochastic} and the references therein for more details.

We see that the stochastic process $X_t$ in Eq.(\ref{p0}) under  L\'{e}vy noise depends on the following parameters. The first one is an uncertain system parameter $\theta$. In general, the estimated parameter $\theta$  plays a key role  in the system model, which could be a bifurcation parameter inducing a transition between states. 
The control parameter $\theta$ could be a greenhouse factor in the case of the energy balance model \cite{kaper2013matclimat}, or a freshwater forcing strength in the thermohaline circulation one\cite{cessi1994simple}. 
Besides, there are uncertain  L\'{e}vy parameters: the non-Gaussian  index $\alpha$ and the L\'{e}vy noise intensity $\epsilon$. 

Let us assume that we have access to a set of observations $y= (y_1, y_2, \cdots, y_n)$, which are the version of the process $X_t$ with a non-Gaussian distribution sampled at discrete times $t_k \in [0, T]$ for $k = 1, 2, \cdots, n$, i.e., $y_k = X_{t_k}$ for $k = 1, 2, \cdots, n$. In this paper, we will discuss the problem of estimating the parameters $\alpha$, $\epsilon$ and $\theta$ simultaneously using the observations $y= (y_1, y_2, \cdots, y_n)$.


\begin{figure}[!t]
\begin{minipage}[b]{0.65 \textwidth}
\leftline{(a)}
\centerline{\includegraphics[height = 5.5cm, width = 11.5cm]{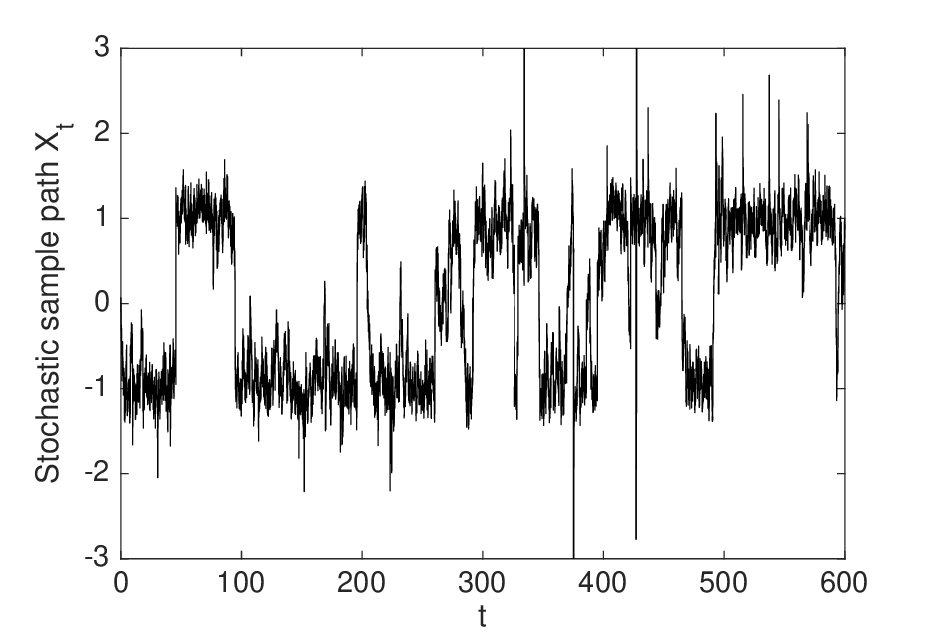}}
\end{minipage}
\hfill
\begin{minipage}[b]{0.35 \textwidth}
\leftline{(b)}
\centerline{\includegraphics[height = 5.5cm, width = 7.5cm]{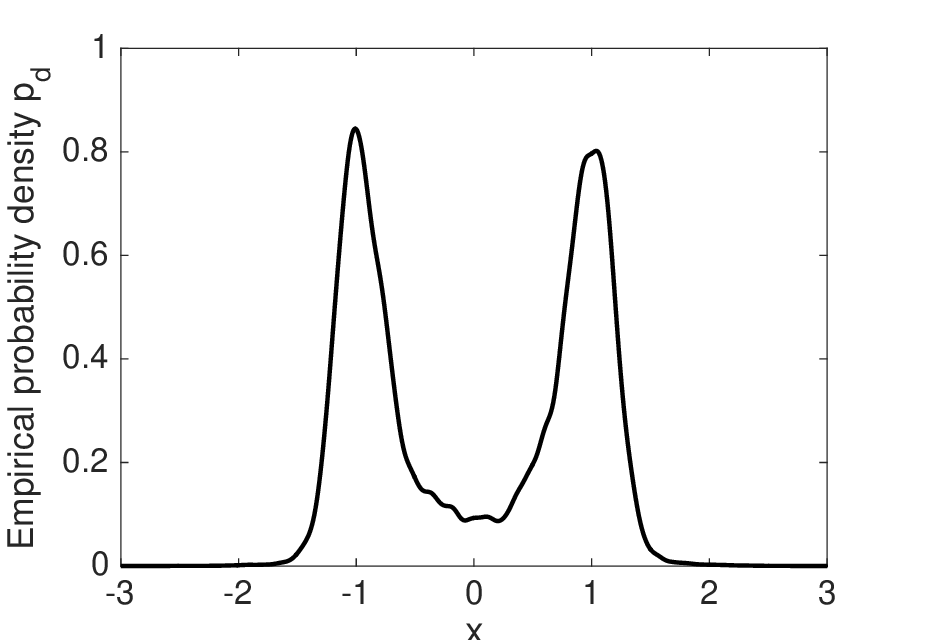}}
\end{minipage}
\caption{(a) The stochastic sample path of the process $X_t$ (Eq.(\ref{sd})) with $\alpha=1.7$, $\epsilon=0.3$ and $\theta = 1$.  (b) The  empirical probability density $p_d$ of the stochastic trajectory (a).}
\label{fig:tra}
\end{figure}

To reach this purpose, we would like to introduce the Hellinger distance. It is used to quantify the similarity between two probability distributions \cite{beran1977minimum}. The Hellinger distance between two probability density functions $p(x)$ and $p_d(x)$ is
\begin{equation*}
H^2(p, p_d)  =  \frac{1}{2}\int_{\mathbb{R}^1}\left(\sqrt{p(x)} - \sqrt{p_d(x)}\right)^2 dx.
\end{equation*}
The Hellinger distance $H$ satisfies the property: $0 \leqslant H(p, p_d) \leqslant 1$. 
Here, $p_d$ is the empirical probability density from an observation data set $y= (y_1, y_2, \cdots, y_n)$.  The probability density function $p(x)$  is a solution of the nonlocal Fokker-Planck equation $p(x, t)$ at time $t$.
\begin{equation}\label{fp1}
\frac{\partial }{\partial t}p(x, t) = -\frac{\partial}{\partial x}(f(x, \theta)p(x, t)) +  \epsilon^{\alpha}\int_{\mathbb{R}^{1}\backslash \{0\}}\left[p(x + y, t) - p(x, t) - I_{|y|<1}(y) \; y \frac{\partial}{\partial x} p(x, t)\right]\nu_{\alpha}(dy). 
  \end{equation}
The integral part in the right hand side is actually the nonlocal Laplacian operator. This nonlocality is the manifestation of effect of non-Gaussian  L\'{e}vy fluctuations \cite{duan2015introduction}. The equation fulfills an initial condition
\begin{equation*}
\lim_{t \rightarrow 0} p(x , t | x_0 , 0) = \delta (x - x_0).
\end{equation*}

We consider that the observation set comes from an $\alpha$-stable L\'{e}vy distribution $p(x, t)$. Associated with each probability density is the parameters set $\lambda = (\theta, \alpha, \epsilon )\in \Theta$, where $\Theta$ is called the parameter space, a finite-dimensional subset of the Euclidean space. Evaluating the Hellinger distance  at the observed data  set $y$  gives  an objective function
\begin{equation*}
G(\lambda) = H^2(p(x, \lambda), p_d(x)).
\end{equation*}

The Hellinger distance estimation aims to find the value of the model parameters that minimize the objective function over the parameter space $\Theta$, that is 
\begin{equation*}
\hat{\lambda} = \arg \min_{\lambda \in \Theta} {G}(\lambda).
\end{equation*}


To address the probability density $p(x, \lambda)$, we use the numerical algorithm of Gao \emph{et al.}\cite{gao2016fokker}  to solve the nonlocal differential equation in Eq.(\ref{fp1}) under the absorbing condition.
This absorbing condition means that the probability of finding ``partical" $X_t$ outside the finite interval  $D = (a, b)$ is zero. We  decompose the integral part of Eq.(\ref{fp1}) into three parts $\int_{\mathbb{R}^1} = \int^{a-x}_{-\infty} + \int^{a-x}_{b-x} + \int^\infty_{b-x}$ in $\mathbb{R}^1$ and analytically evaluate the first and third integrals, then Eq.(\ref{fp1}) changes to
\begin{align}\label{p3}
\frac{\partial }{\partial t}p(x, t) &=  -\frac{\partial }{\partial x}(f(x)p(x, t)) - \frac{\epsilon^\alpha C_\alpha}{\alpha} \left[\frac{1}{(x-a)^\alpha} + \frac{1}{(b-x)^\alpha}\right] p(x, t)     \notag\\
&+ \epsilon^\alpha C_\alpha \int_{a-x}^{b-x} \frac{p(x+y, t) - p(x, t)- I_{|y|<1}(y)y \frac{\partial}{\partial x}p(x, t) }{|y|^{1 + \alpha}}dy,
 \end{align}
for $x \in (a, b)$. 
 The non-Gaussian  index  $\alpha \in (0, 2)$ and the L\'evy intensity  $\epsilon \in (0, 1]$. 
 
 \begin{figure}[!t]
\begin{minipage}[b]{0.5 \textwidth}
\leftline{(a)}
\centerline{\includegraphics[height = 6cm, width = 9.4cm]{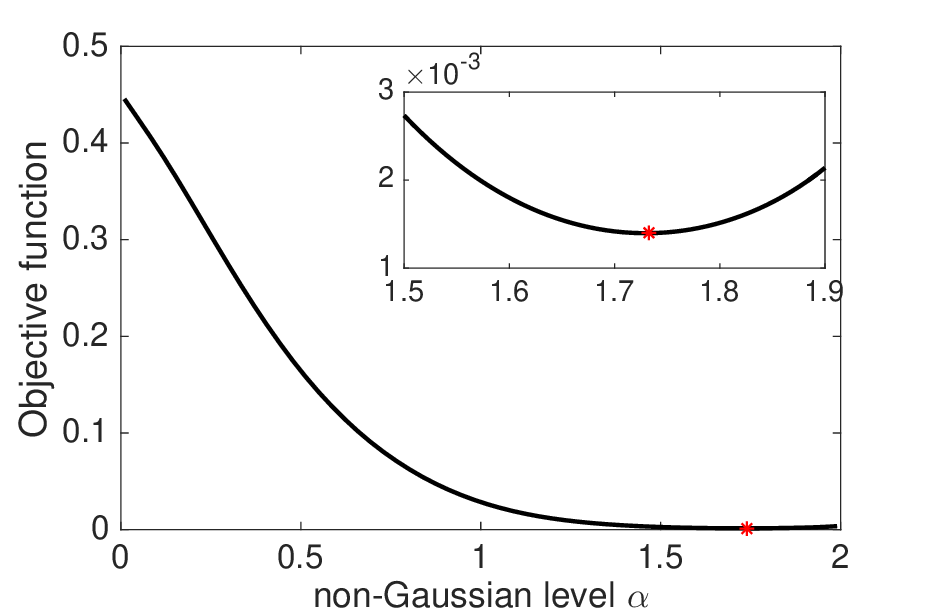}}
\end{minipage}
\hfill
\begin{minipage}[b]{0.5 \textwidth}
\leftline{(b)}
\centerline{\includegraphics[height = 6cm, width = 9.4cm]{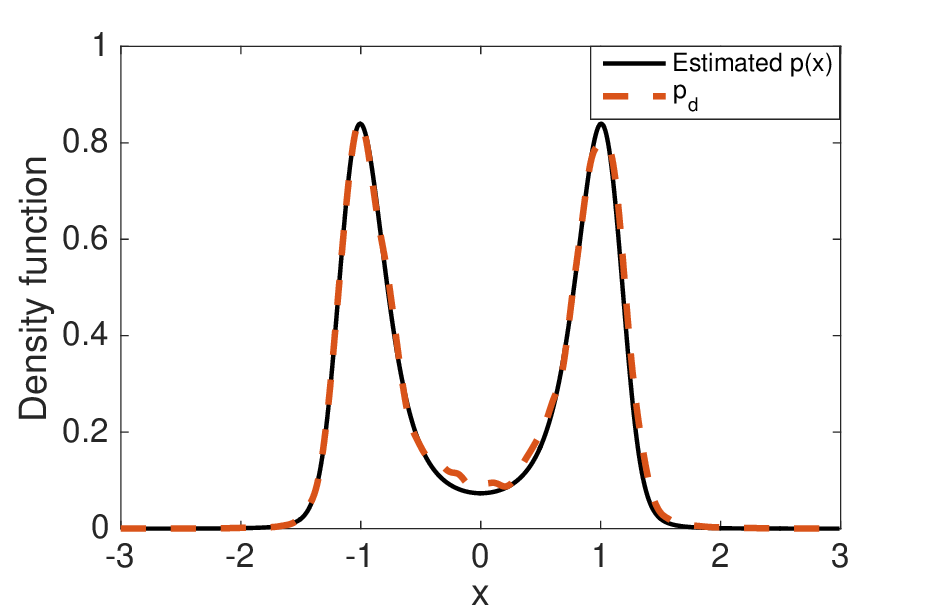}}
\end{minipage}
\caption{(a) The Hellinger distance estimation of $\alpha$. The estimated value $\hat{\alpha} = 1.738$ corresponds to the  minimum objective function $G _\alpha= 0.0014$ (red $\star$).  (b) The estimated density function $p(x)$  with $\hat{\alpha} = 1.738$, $\theta = 1$ and $\epsilon = 0.3$ compared with the empirical density  $p_d$ in Fig.\ref{fig:tra}(b).}
\label{fig:alp}
\end{figure}

\section{Numerical Experiments}
We now explore how the Hellinger distance can be used to estimate the parameters of SDE driven by a symmetric $\alpha$-stable L\'{e}vy process. We consider the following example
\begin{equation}\label{sd}
{\rm d}X(t) = (-\theta X(t)^3 + X(t)){\rm d}t+ \epsilon {\rm d}L^\alpha(t) , \quad X(0)=0 \in \mathbb{R}^1,
\end{equation}
In this example, a nonlinear drift term is  $f(x, \theta) = -\theta x^3 + x$ with uncertain parameter $\theta$. We start with training data from numerical simulations of Eq.(\ref{sd}).  The  stochastic trajectory can be regarded as a heavy-tailed time series with the parameters $\alpha=1.7$, $\theta=1$ and $\epsilon=0.3$   shown in the Fig. \ref{fig:tra}(a). \textcolor[rgb]{0,0,0}{Here, the choice of parameters $\alpha \in (0, 2)$, $\theta \in [0.5, 1.5]$ and $\epsilon \in (0, 1]$  are arbitrary. Our work focus on the comparison the similarity between the estimated  and original parameters in the SDE (\ref{sd}) by minimizing the Hellinger distance. The Hellinger distance measures  the similarity between an empirical probability density of non-Gaussian observations and a solution (as a probability density) of the associated nonlocal Fokker-Planck equation (\ref{fp1}).}

The empirical probability density function $p_d$ for observations  (Fig. \ref{fig:tra}(a)) could be determined by the normal kernel method. In the simulation, we use the MATLAB function \emph{ksensity} to evaluate the $p_d(x)$ for $x \in (-3, 3)$ as shown in Fig. \ref{fig:tra}(b). A selected bandwidth is $h=1.8s/n^{1/5}$, where $n = 6\times10^6$ is the number of observed data points and $s$ is the standard deviation of the data set.

 Now we provide the details of the computation to infer parameters from the observations. The value of the  parameters are estimated by minimizing the Hellinger distance  over the parameter space $\Theta$. We shall first estimate single parameter assuming that the values of the other parameters are known, and then estimate multiple parameters 
  
 \subsection{Estimation for a single parameter}\label{s31}
We want to find out an estimation of $\alpha$ by achieving a numerical optimization of the objective function $G$  of the Hellinger distance.  We consider the parameter  $\lambda =\alpha $ is in the parameter space $\Theta =(0, 2) \subset  \mathbb{R}^1$,  and assume that the other parameters are known, i.e. $\epsilon = 0.3$ and $\theta = 1$. Then the objective function $G$  of the Hellinger distance is
\begin{equation*}
G(\alpha)  =  \frac{1}{2}\int_{\mathbb{R}^1}\left(\sqrt{p(x, \alpha)} - \sqrt{p_d(x)}\right)^2 dx.
\end{equation*}
Based on the numerical algorithm of Gao \emph{et al.}\cite{gao2016fokker},  the probability density $p(x, \alpha)$ is  solved by the nonlocal differential equation (\ref{fp1}) given  $\alpha \in (0, 2) $ for $x \in D= (-3, 3)$ at $t=50$.  In the numerical simulations,  the probability profile of its initial position is   Gaussian  $p(x,0) = \sqrt{\frac{40}{\pi}}e^{-40x^2}$. We have chosen the spatial resolution $h=0.003$ and the time step size $\Delta t = 0.5h^2$. 


In Fig. \ref{fig:alp}(a),  we employ a discretization step of $\Delta \alpha = 0.035$ and use 55 grid points for $\alpha \in (0, 2)$. Then the estimation of  $\hat{\alpha} = 1.738$ is obtained with the minimum value of the Hellinger distance $G (\hat{\alpha})= 0.001399$ over  the parameter space $\Theta = (0, 2)$. Furthermore, we restrict $\alpha$ on a small region $ [1.5, 1.9]$ for accurately estimation.
The  result illustrates that the  estimated values of  $\hat{\alpha} \in [1.684, 1.764]$ contains the true value of $\alpha=1.7$ with the Hellinger distance $G(\hat{\alpha}) < 0.00145$ (inset figure in Fig. \ref{fig:alp}(a)) .  It means that we could find an optimal interval for the estimated parameter $\alpha$.  As an illustration, we show  the results of the probability density $p(x, \alpha)$ of SDE (\ref{sd}) with estimated value $\hat{\alpha} = 1.738$ (dashed)  and  the empirical density $p_d$ from the observed data (dotted) in Fig. \ref{fig:tra}(b). We can see that the estimated probability density presents the goodness-of-fit to the empirical one.

 \begin{figure}[t]
\begin{minipage}[b]{0.5 \textwidth}
\leftline{(a)}
\centerline{\includegraphics[height = 6cm, width = 9.4cm]{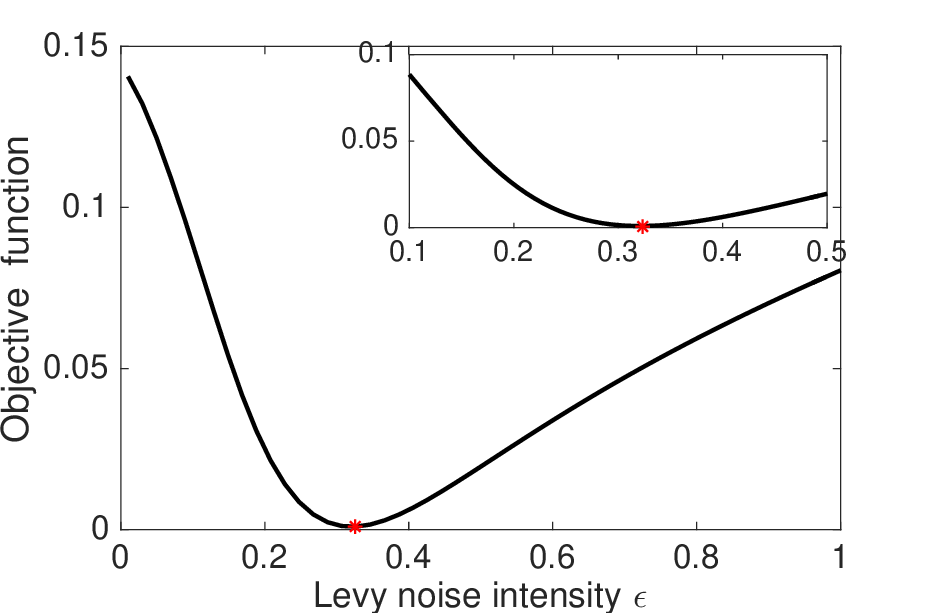}}
\end{minipage}
\hfill
\begin{minipage}[b]{0.5 \textwidth}
\leftline{(b)}
\centerline{\includegraphics[height = 6cm, width = 9.4cm]{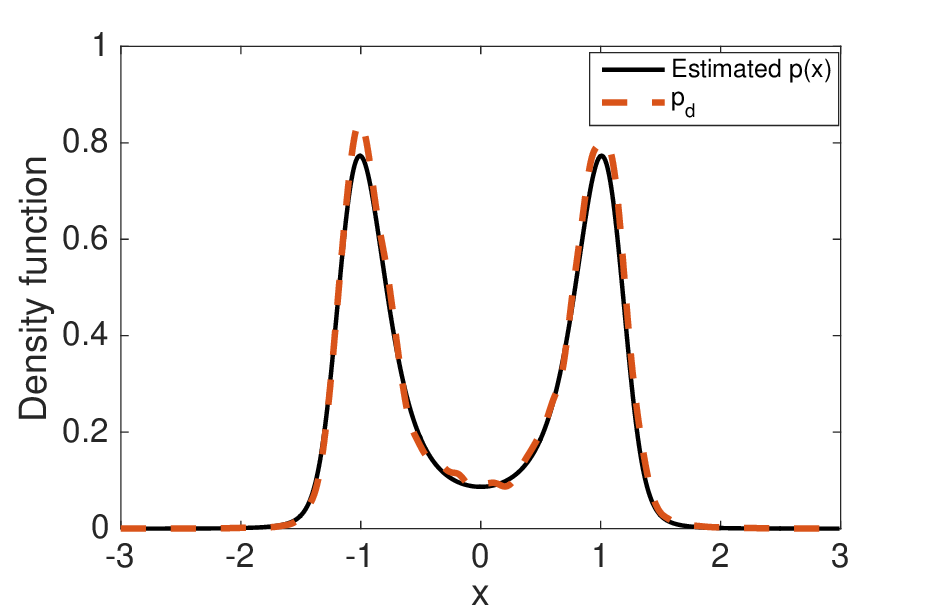}}
\end{minipage}
\caption{(a) The Hellinger distance estimation of $\epsilon$. The estimated value $\hat{\epsilon} = 0.3268$ corresponds to the  minimum objective function $G (\hat{\epsilon}) = 0.000951$ (red $\star$).  (b) The estimated density function $p(x)$  with $\hat{\epsilon} = 0.3268$, $\theta = 1$ and $\alpha = 1.7$ compared with the empirical density  $p_d$ in Fig.\ref{fig:tra}(b).
}
\label{fig:eps}
\end{figure}

Similarity, we explore the dependence of the objective function $G$ on the value of $\epsilon$  keeping the other parameters  fixed. 
Fig. \ref{fig:eps}(a) shows the minimized the  $G (\hat{\epsilon}) = 0.0009$ of the Hellinger distance that corresponds to the estimated value  of the L\'{e}vy noise intensity $\hat{\epsilon }= 0.3268$.   
 Meanwhile, we could get the optimal estimation interval $\hat{\epsilon} \in [0.292, 0.348]$ for the Hellinger distance $G(\hat{\epsilon} )< 0.0019$. This domain  includes the true value of $\epsilon = 0.3$.
In Figure \ref{fig:eps}(b), the empirical density $p_d$ is well fitted by  the probability density $p(x, \epsilon)$ with $\hat{\epsilon }= 0.3268$.

\textcolor[rgb]{0,0,0}{We have inferred the parameters $\alpha$ and $\epsilon$ by considering the Hellinger distance, respectively. Next, we would like to compare the Hellinger distance with other commonly used metrics, such as the $L^2$ norm distance, the maximum absolute error  distance and the S\o rensen distance to quantify the similarity between two probability distributions. The  objective function of the $L^2$ norm distance is defined as
\begin{equation*}
G(\lambda) = \frac{\|p(\lambda, x) - p_{d}(x)\|^2_2}{\|p_{d}(x)\|^2_2}.
\end{equation*}
 While the  maximum absolute approximation distance is
\begin{equation*}
G(\lambda) =  \max| p(\lambda, x) - p_{d}(x)|.
\end{equation*}
The S\o rensen distance is used in ecology model. The expression of the objective function for the S\o rensen distance is
\begin{equation*}
G(\lambda) =  \frac{\| p(\lambda, x) - p_{d}(x)\|_1 }{\| p(\lambda, x) + p_{d}(x)\|_1}.
\end{equation*}
Here the  $\| z\|_n = (|z_1|^n+ |z_2|^n + \cdots + |z_m|^n)$ is  defined as $n$-norm of $z = (z_1, z_2, \cdots, z_m)$. The estimation of  uncertain parameters set $\lambda$ could be achieved by minimizing $G$, i.e.,  $\hat{\lambda} = \arg \min G(\lambda)$ for $\lambda \in \Theta$.}


\textcolor[rgb]{0,0,0}{Next, let us examine the effect of these kinds of distances on the estimation of $\alpha$ and $\epsilon$, respectively. We keep the other parameters and the divided subintervals   the same as those in Figs. \ref{fig:alp} and \ref{fig:eps}.
In Table \ref{T1}, the results on these distance show that the Hellinger distance gives a better estimation for $\alpha$ than the others. In contrast,  all four distances show a good fit to the true  $\epsilon$. In this example, the Hellinger distance is the most effective method to estimate parameters.}

\begin{table}[htp]
\begin{center}
\begin{tabular}{ccccccc}
\hline
Distance & True value $\alpha$ & Estimated  $\hat{\alpha}$ & True value $\epsilon$ & Estimated $\hat{\epsilon}$ & $G_{\alpha}$ &  $G_{\epsilon}$\\
\hline
Hellinger   &1.7& 1.7380 & 0.3 & 0.3268 & 0.0014 & 0.0009\\
$L^2$ norm  &1.7& 1.8100 & 0.3 & 0.3070 & 0.0028 & 0.0032\\
Maximum absolute error   &1.7& 1.9180 & 0.3 & 0.3070 & 0.0427 & 0.062\\
S\o rensen  & 1.7 & 1.7740 & 0.3 &  0.3070 &  0.0299 & 0.0295\\

\hline
\end{tabular}
\end{center}
\caption{\textcolor[rgb]{0,0,0}{Compared with different distances: Hellinger distance, $L^2$ norm distance, maximum absolute error distance and the S\o rensen distance.}}
\label{T1}
\end{table}%

 \begin{figure}[!t]
\begin{minipage}[b]{0.5 \textwidth}
\leftline{(a)}
\centerline{\includegraphics[height = 6cm, width = 9.4cm]{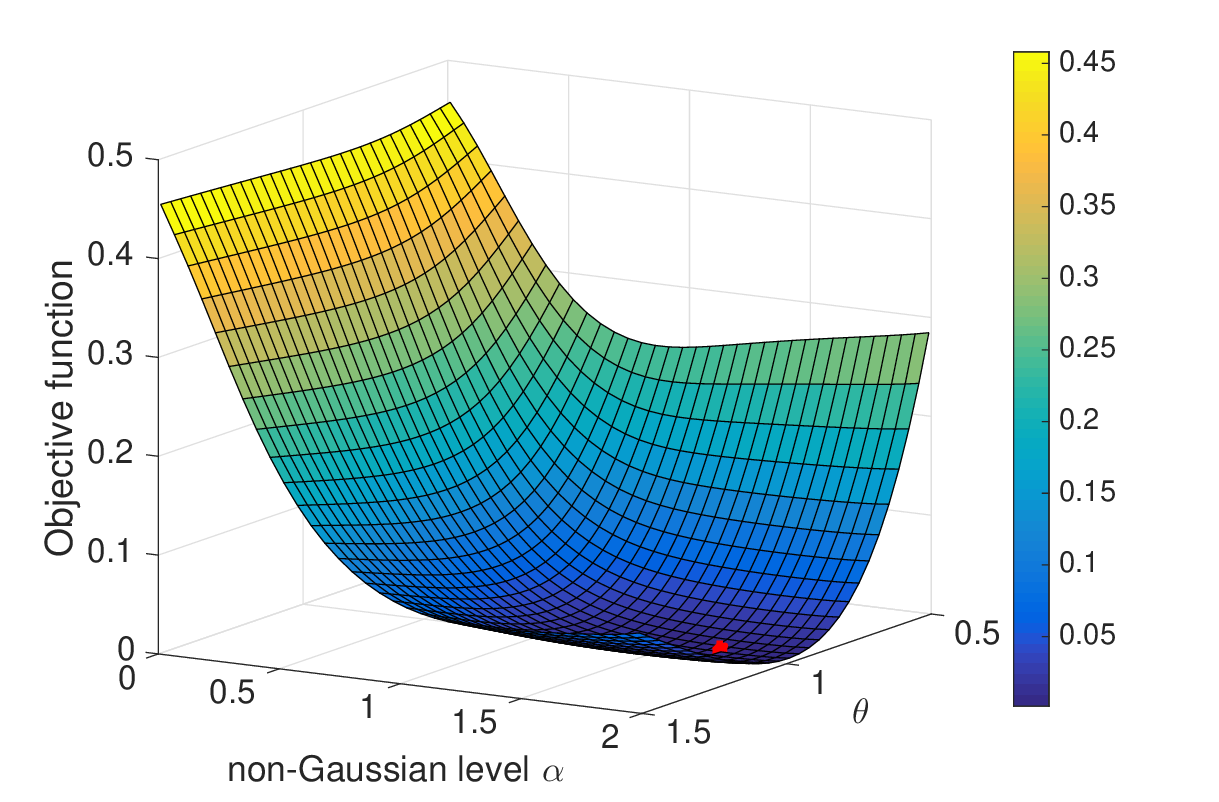}}
\end{minipage}
\hfill
\begin{minipage}[b]{0.5 \textwidth}
\leftline{(b)}
\centerline{\includegraphics[height = 6cm, width = 8.5cm]{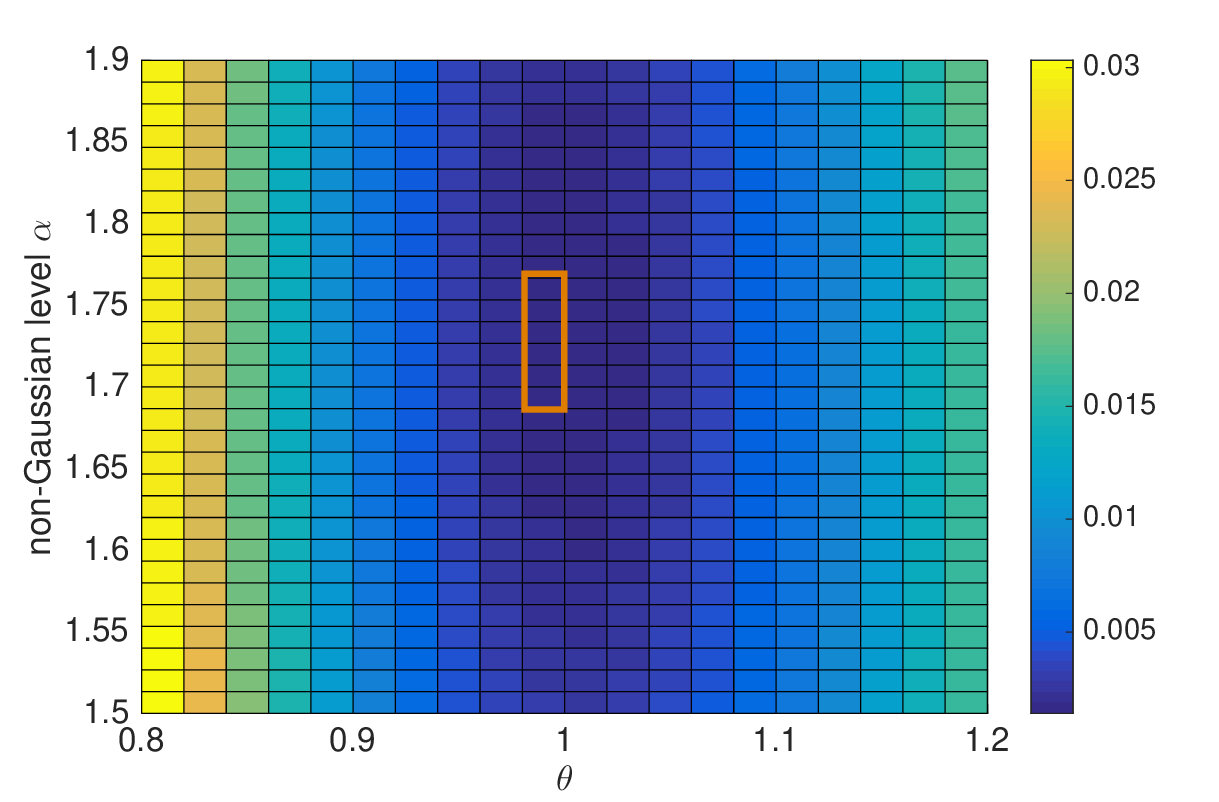}}
\end{minipage}
\caption{(a) Estimation of  $\alpha$ and $\theta$   by optimization of the Hellinger distance.  The estimated values are $\hat{\alpha} = 1.707$ and   $\hat{\theta} = 0.9828$ with the  minimum $G(\hat{\alpha}, \hat{\theta})= 0.0014$ (red $\star$). (b)  The optimal estimation domain of $\hat{\alpha}$ and $\hat{\theta}$ is $  [1.6867, 1.78] \times [0.98, 1]$.}
\label{fig:alpha}
\end{figure}

 \subsection{Estimation for multiple parameters}

The above example has verified that our method is feasible for estimating a single parameter by minimizing the Hellinger distance. Next, we will apply this approach to  estimate the multiple unknown parameters.  First, we simplify our model by assuming that one parameter $\epsilon = 0.3$  is known and then estimate  $\alpha$ and $\theta$, while keeping the other factors the same as in the  section \ref{s31}. The objective function $G$ is given by Hellinger distance
\begin{equation*}
G(\alpha, \theta)  =  \frac{1}{2}\int_{\mathbb{R}^1}\left(\sqrt{p(x, \alpha, \theta)} - \sqrt{p_d(x)}\right)^2 dx.
\end{equation*}
The probability density $p(x, \alpha, \theta)$ is a solution of the nonlocal differential equation (\ref{fp1}) given values of $\alpha \in (0, 2)$ and $\theta \in [0.5, 1.5]$  at $t=50$. 
Fig. \ref{fig:alpha}(a) shows that the objective function $G$  changes with the  values of $\alpha$ and $\theta$ in the parameter space $\Theta =(0, 2) \times  [0.5, 1.5] \subset  \mathbb{R}^2$. The  minimum value of  $G(\hat{\alpha}, \hat{\theta})= 0.0014$ is identified with $\hat{\alpha} = 1.707$ and   $\hat{\theta} = 0.9828$. In the same manner,   the optimal estimation domain of $\hat{\alpha}$ and $\hat{\theta}$ is $ [1.6867, 1.78]\times [0.98, 1]$ (orange rectangular frame) by further restricting the  range of parameters as shown in Fig. \ref{fig:alpha}(b).

 \begin{figure}[!t]
\begin{minipage}[b]{0.5 \textwidth}
\leftline{(a)}
\centerline{\includegraphics[height = 6cm, width = 9.4cm]{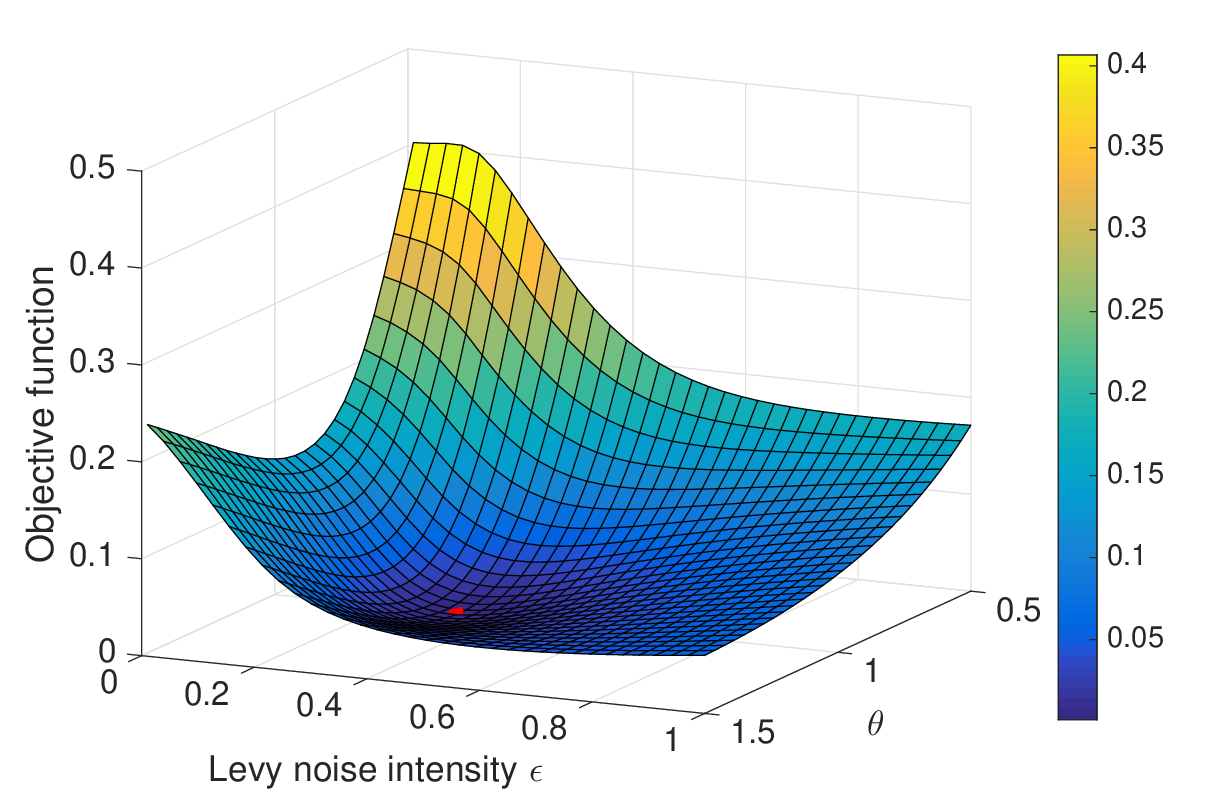}}
\end{minipage}
\hfill
\begin{minipage}[b]{0.5 \textwidth}
\leftline{(b)}
\centerline{\includegraphics[height = 6cm, width = 8.5cm]{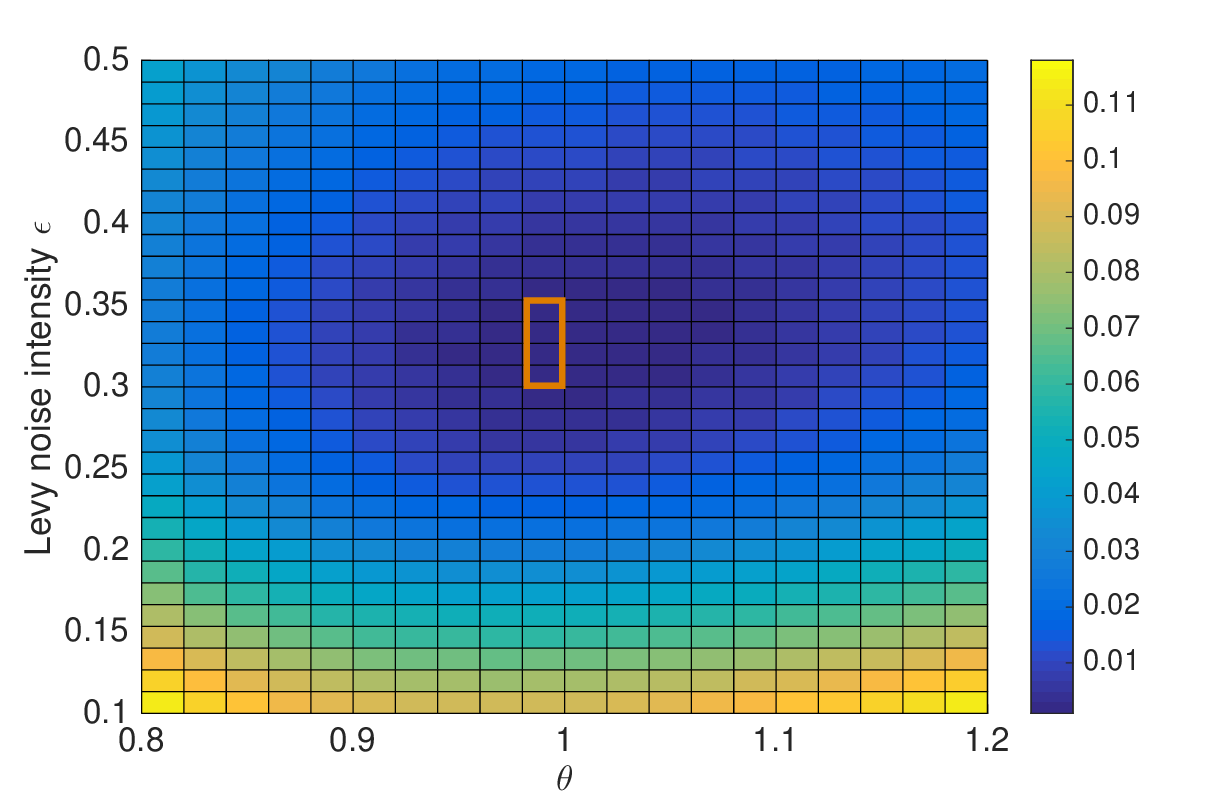}}
\end{minipage}
\caption{(a) Estimation of  $\epsilon$  and $\theta$.  The estimated values are $\hat{\epsilon}= 0.3303$ and $\hat{\theta} = 1.0185$ with  the  minimum $G(\hat{\epsilon}, \hat{\theta})= 0.0012$ (red $\star$). (b) The optimal estimation domain is $ [0.3, 0.3267] \times [0.9, 1.02]$.}
\label{fig:epsa}
\end{figure}

 \begin{figure}[!t]
\begin{minipage}[b]{0.5 \textwidth}
\leftline{(a)}
\centerline{\includegraphics[height = 6cm, width = 9.4cm]{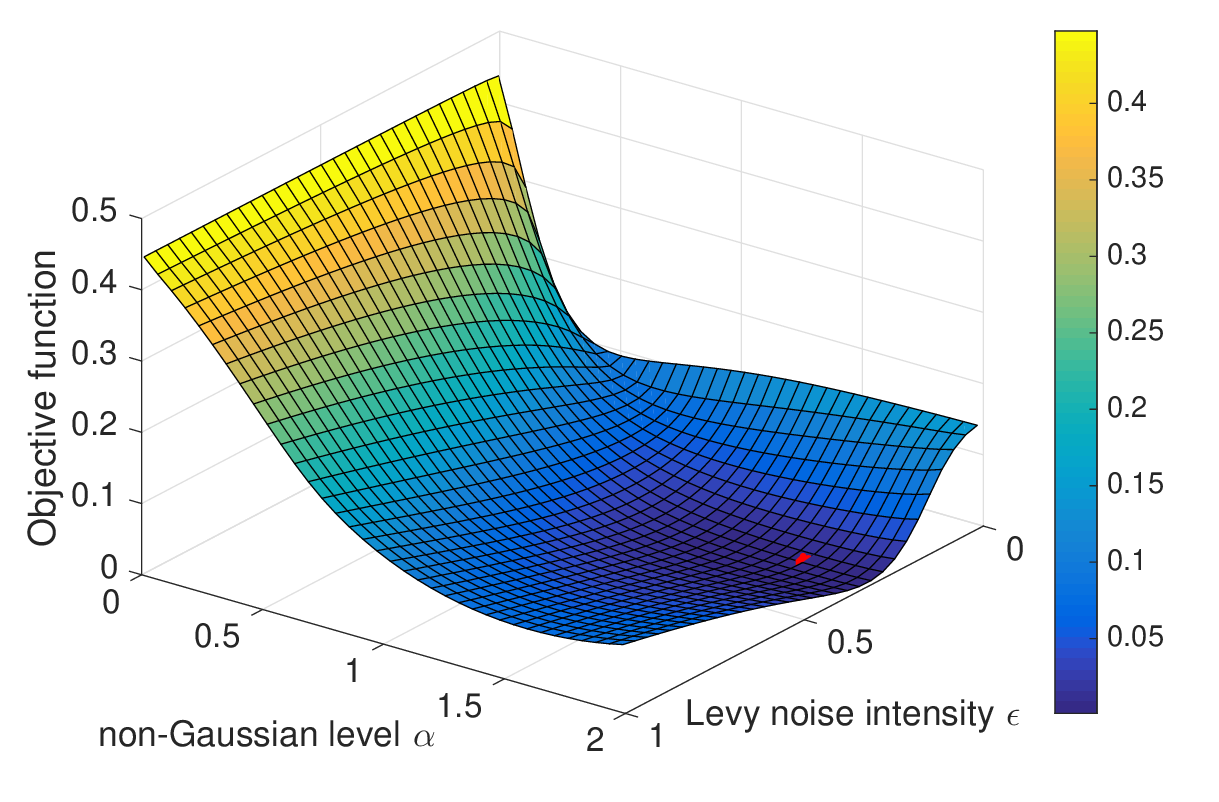}}
\end{minipage}
\hfill
\begin{minipage}[b]{0.5 \textwidth}
\leftline{(b)}
\centerline{\includegraphics[height = 6cm, width = 8.5cm]{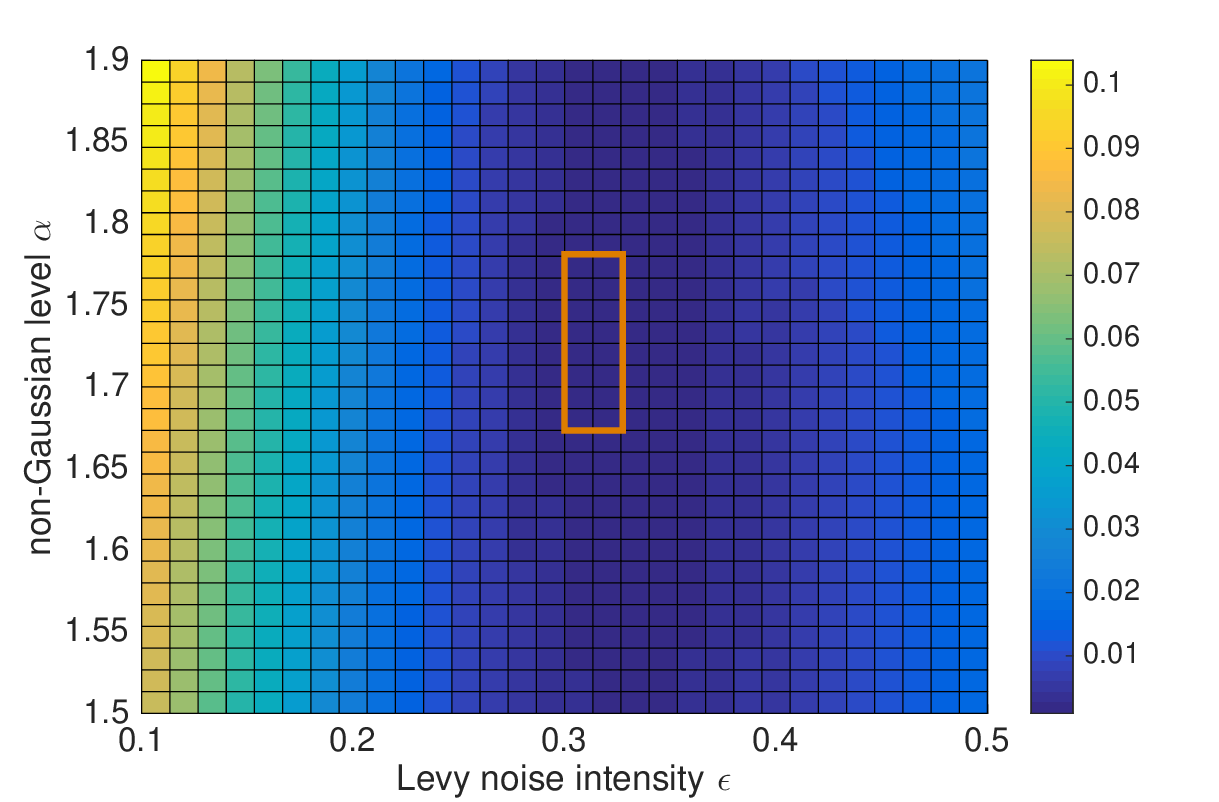}}
\end{minipage}
\caption{(a)  Estimation of $\alpha$ and $\epsilon$.   The estimated values are $\hat{\alpha}= 1.707$ and   $\hat{\epsilon}= 0.307$ with  the  minimum $G(\hat{\alpha}, \hat{\epsilon})= 0.0011$ (red $\star$). (b) The optimal estimation domain is  $ [1.6733, 1.78] \times  [0.3133, 0.3267]$.}
\label{fig:alpeps}
\end{figure}

\begin{table}[htp]
\begin{center}
\begin{tabular}{ccccc}
\hline
Parameter & True value & Estimated $\hat{\lambda}_1 = ( \hat{\epsilon}, \hat{\theta} )$ & Estimated $\hat{\lambda}_2= (\hat{\alpha}, \hat{\epsilon} )$\\
\hline
$\theta$  & 1.0 &  1.0185  & --\\
$\alpha$ & 1.7 &   --  & 1.707\\
$\epsilon$ & 0.3 & 0.3303 & 0.307 \\
\hline
$G$ & 0 & 0.0012 & 0.0011\\
\hline
\end{tabular}
\end{center}
\caption{Estimation of $\lambda_1 = (\epsilon, \theta)$ and $\lambda_2= (\alpha, \epsilon )$.}
\label{T2}
\end{table}%

Second, we take into account the estimation of the other two combinations of all three parameters, $\lambda_1 = (\epsilon, \theta )$ and $\lambda_2= (\alpha, \epsilon )$ corresponding to the  parameter spaces $\Theta_1 =  (0, 1]\times [0.5, 1.5] $,  $\Theta_2 =  (0, 2) \times  (0, 1]$, respectively.
The estimated results  are found  by  the minimized  the Hellinger distance as shown in Table \ref{T2}.  Meanwhile, we can also determine the optimal domains of the estimated parameters sets $\hat{\lambda}_1$ and $\hat{\lambda}_2$ as shown in Figs. \ref{fig:epsa}(b) and \ref{fig:alpeps}(b). The result shows that the Hellinger distance $G(\hat{\alpha}, \hat{\epsilon}) \leqslant 0.001$ if the estimated parameters set $\hat{\lambda}_1$ belongs to the domain $[1.6733, 1.78] \times [0.3, 0.3267]$. Meanwhile, the estimation domain of $\hat{\lambda}_2$ is  $[0.3, 0.34] \times [0.98, 1]$ if the Hellinger distance $G(\hat{\epsilon}, \hat{\theta}) \leqslant 0.0015$.

\begin{figure}[h]
\centerline{\includegraphics[width=14cm ,height= 10cm]{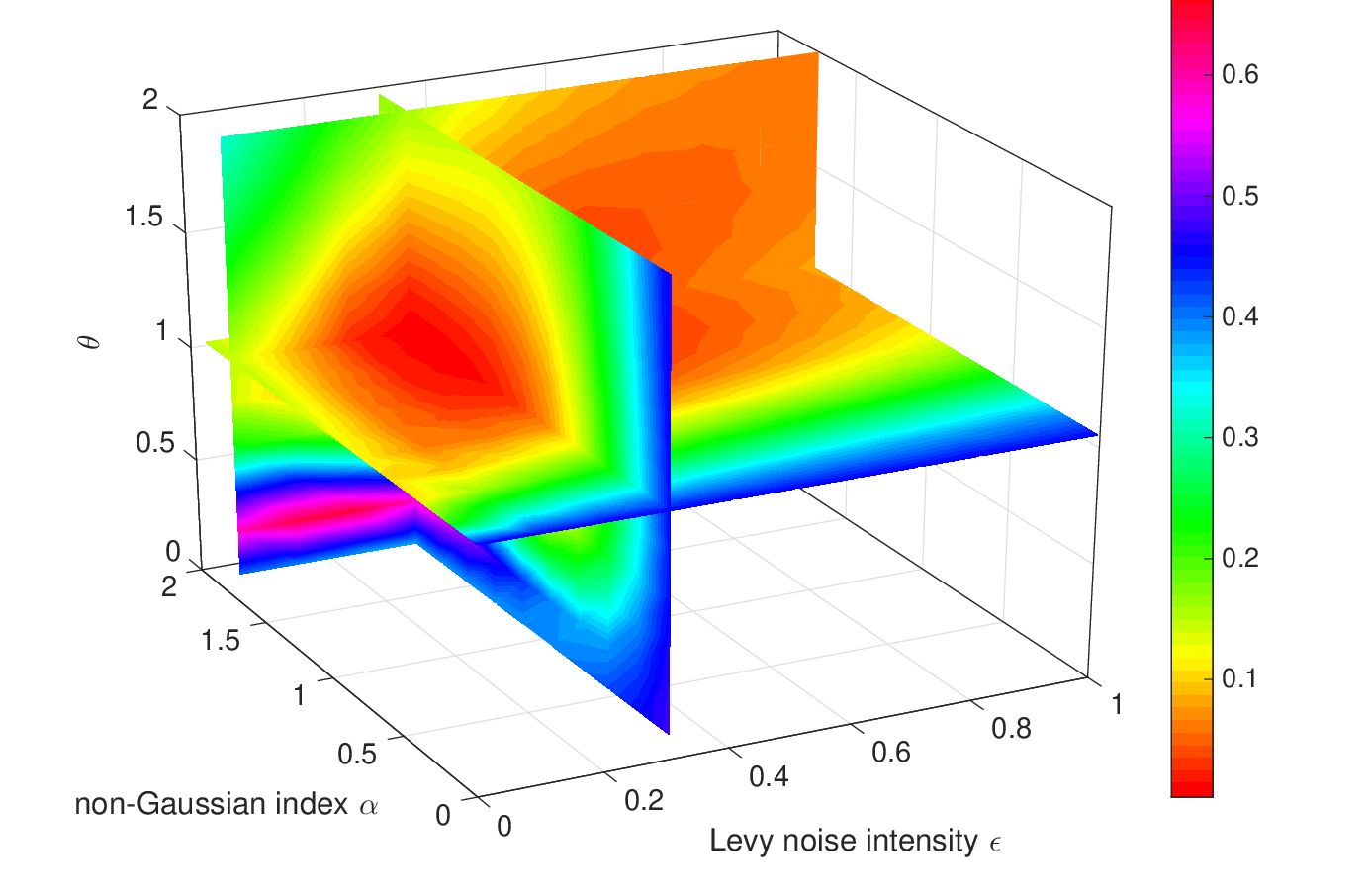}}
\caption{ The Hellinger distance estimation  of $\alpha$, $\epsilon$  and $\theta$.  The estimated $\hat{\epsilon} = 0.307$, $\hat{\theta} = 1.05$ and $\hat{\alpha} = 1.7380$ for the  minimum $G(\hat{\alpha}, \hat{\theta}, \hat{\epsilon})= 0.0027$.}
\label{fig:all}
\end{figure}

Finally, we seek all these parameters $\lambda= (\alpha, \epsilon, \theta )$ such that the Hellinger distance reaches the minimum value in the  parameter spaces $\Theta =  (0, 1] \times (0, 2] \times (0, 2) \subset \mathbb{R}^3$. 
\begin{equation*}
G(\theta,\epsilon, \alpha)  =  \frac{1}{2}\int_{\mathbb{R}^1}\left(\sqrt{p(x, \theta,\epsilon, \alpha)} - \sqrt{p_d(x)}\right)^2 dx.
\end{equation*}
The  values of the model parameters $\hat{\epsilon} = 0.307$, $\hat{\theta} = 1.05$ and $\hat{\alpha} = 1.7380$ are achieved  by minimizing the Hellinger distance $G(\hat{\alpha}, \hat{\theta}, \hat{\epsilon})= 0.0027$ over the parameter space $\Theta \subset\mathbb{R}^3$.  Estimated results  $\hat{\epsilon} = 0.307$, $\hat{\theta} = 1.05$ and $\hat{\alpha} = 1.7380$ defines a slice plane in the $\epsilon$-axis, $\theta$-axis, or $\alpha$-axis direction as shown in Fig. \ref{fig:all}.

\section{Conclusion}

In summary, we consider a non-Gaussian dynamical system containing uncertain parameters. An approach of parameter estimation is proposed by  numerical optimization of the Hellinger distance between two probability distributions. The one probability density  $p(x)$  is a solution of the nonlocal Fokker-Planck equation at time $t$ for a stochastic dynamical system $X(t)$ driven by an $\alpha$-stable L\'{e}vy process. The other one is  the empirical probability density $p_d$  from observations data of discrete version of the process.   The approach is used to find all out the estimation  of  single parameter and multiple  parameters, by a numerical optimization of  the Hellinger distance over the parameters space.  The results of an example  verified that this method is feasible for estimating non-Gaussian parameters $\alpha$, $\epsilon$ and other system parameter by the Hellinger distance.  Compared with the $L^2$ norm, the maximum absolute error distances and S\o rensen distance, the Hellinger distance is the most effective method to estimate parameters in this example. Meanwhile, we could find an optimal interval for the estimated parameters.

This  approach can be used to establish parameter estimations for a data-driven dynamical system, the observations data with jumps and heavy-tailed distribution. A very important future work will be a model study of the abrupt climate changes in the Dansgaard-Oeschger events  with non-Gaussian distribution. The approach would be applied  to estimate the system parameters and non-Gaussian parameters in this model.

\textcolor[rgb]{0,0,0}{\section*{Data Availability}
All computational results are implemented with MATLAB R2015b running in an Intel Xeon(R)  CPU E5-2667 v4 @ 3.20 GHz machine. The data that support the findings of this study are available in GitHub at: https://github.com/yayun55/Estimating-uncertainty-for-the-observed-non-Gaussian-data.}


\section*{Declaration of Competing Interest}
The authors declare that they have no known competing financial interests or personal relationships that could have appeared to influence the work reported in this paper. 

\section*{Author Contributions} Y. Zheng designed the research. Y. Zheng and  F. Yang  performed computations and wrote the first draft of the manuscript. J. Kurths and J. Duan analysed the results and concepts development. All authors conducted research discussions and reviewed the manuscript.

\section*{Acknowledgements}
We would like to thank Xiaoli Chen, Xiujun Cheng and Yang Liu for discussions about computation. This work is supported by the
National Natural Science Foundation of China (grants No.11801192 and No. 11771449), Jiangsu University Project Grant (No. 20JDG071) and Russian Ministry of Science and Education "Digital biodesign and personalised healthcare”.


\end{document}